\documentstyle[11pt]{article}
\newtheorem{tm}{Theorem}[subsection]
\newtheorem{lm}[tm]{Lemma}
\newtheorem{pr}[tm]{Proposition}
\newtheorem{rmk}[tm]{Remark}
\newtheorem{cor}[tm]{Corollary}

\newtheorem{??}[tm]{Question}
\newtheorem{defi}[tm]{Definition}

\font\tenmsb=msbm10
\font\sevenmsb=msbm7
\font\fivemsb=msbm5

\newfam\msbfam
\textfont\msbfam=\tenmsb
\scriptfont\msbfam=\sevenmsb
\scriptscriptfont\msbfam=\fivemsb
\def\Bbb#1{{\fam\msbfam #1}}

\font\teneufm=eufm10
\font\seveneufm=eufm7
\font\fiveeufm=eufm5
\newfam\eufmfam
\textfont\eufmfam=\teneufm
\scriptfont\eufmfam=\seveneufm
\scriptscriptfont\eufmfam=\fiveeufm

\def\lorw{\longrightarrow}

\newcommand\ci{\cite}

\newcommand\rat{{\Bbb Q}}
\newcommand\comp{{\Bbb C}}

\newcommand\pn[1]{{\Bbb P}^{#1}}

\newcommand\blacksquare{{\hspace*{\fill} $\Box$}}

\title{The Hard Lefschetz Theorem and the topology of semismall maps}
\author{
Mark Andrea A.  de Cataldo\thanks{
Partially supported by N.S.F. Grant DMS 9701779.}\, 
and Luca Migliorini\thanks{ member of GNSAGA, 
supported by MURST funds}
}
\date{\today}

\begin{document}\maketitle

\begin{abstract}
We introduce the notion of {\em lef}$\,$ line bundles on a complex projective 
manifold. We prove that lef line bundles 
satisfy the Hard Lefschetz Theorem, the Lefschetz Decomposition
and the Hodge-Riemann Bilinear Relations.
We study proper holomorphic semismall
maps from complex manifolds and prove that, for constant
coefficients, the Decomposition Theorem is equivalent
to the non-degeneracy of certain intersection forms.
We give a  proof of  the Decomposition Theorem 
for the complex direct image of the constant sheaf when 
the domain and the target are   projective
by proving that
the forms in question  are non-degenerate. A new feature
uncovered by our proof is that the forms are definite.
\end{abstract}

\tableofcontents

\section{Introduction}
Let $X$ be a complex projective $n$-fold endowed with a line bundle
$L$.
The classical statement of the Hard Lefschetz Theorem
is that if $L$ is ample then, for every $r\geq 0$, 
cupping with $c_1( L)^r$ gives an isomorphism 
of $H^{n-r}(X, \rat)$ onto $H^{n+r}(X,\rat)$.

\smallskip
In this paper we introduce the notion of {\em lef} line bundles
on projective manifolds.
They are precisely the line bundles with the property that
$L^{\otimes t}= f^*A$, where $f: X \to Y$ is a projective semismall 
morphism, $t$ is some positive integer and $A$ is an ample line bundle.
Recall that $f$ is {\em semismall} when 
the set  
of points in $Y$, whose fiber has  dimension $k$, has  codimension at least $2k$. 
One has the following strict implications: ample $\Rightarrow$ lef
$\Rightarrow$ semiample$+$big; the results we prove for lef line bundles
fail for the last class.

\smallskip
Our first result is Theorem 
\ref{mainthm}: if $L$ is a line bundle such that
a positive power is generated by its global sections, then
the Hard Lefschetz Theorem holds for $L$ {\em if and only if} $L$ is lef.
See also Proposition \ref{nohlnoss}.
A Lefschetz-type decomposition for $H^*(X, \rat )$ ensues
and we prove the Hodge-Riemann Bilinear relations,  
i.e. the signature of the intersection form
on the associated primitive spaces has 
properties completely analogous to the ample case.
 We offer two proofs.
One is by induction on hyperlane sections; see $\S$\ref{firstproof}.
Another one is by the use of pencils; see $\S$\ref{proofhl}. While the first
is simpler, the second one characterizes monodromy invariants.

\smallskip
A second result  is Theorem \ref{fai} where  we prove
 a higher dimensional generalization 
of the Grauert criterion for the contractibility of a finite set
of curves on a surface: if 
$f:X \to Y$ is a semismall map between $2m$-dimensional projective varieties,
$X$ is nonsingular, 
$y   \in Y  $ is a point and $\dim { f^{-1} (y  ) } =m$, then the intersection form 
associated with the irreducible components of maximal dimension 
of $f^{-1}(y  )$ 
is $(-1)^m$-definite. We call this property the {\em Hodge Index Theorem
for semismall maps}. It is proved using Theorem \ref{mainthm}, the characterization of monodromy invariants
 and  Goresky-MacPherson
version of the Weak Lefschetz Theorem.

\medskip
This result has a surprising  consequence when applied 
  to study  the topology of semismall maps.

\noindent
Beilinson-Bernstein-Deligne, \ci{bbd}, proved 
the following fundamental result, called
the Decomposition Theorem:
if   $g: U \to V$ is a proper morphism of algebraic varieties, then
the direct image complex of an intersection cohomology complex
of geometric origin on $U$, e.g. the intersection cohomology complex 
of $U$, 
 is isomorphic, in the derived category, to a direct sum
of {\em shifted} intersection cohomology complexes. An analogous
result holds over any field. Even over the complex numbers,
finite fields are necessary to their proof. 
Saito, \ci{sa}, proved and extended these results in a K\"ahler
setting using mixed Hodge modules.
 Borho-MacPherson,
\ci{b-m},
have observed that, in the statement of the Decomposition
Theorem,  the shifts disappear when  $g$ is semismall, 
$U$ is smooth and constant coefficients are taken.

\smallskip
Let $f: X \to Y$ be a proper {\em semismall} holomorphic map from a complex manifold.
We propose a new approach
to the Decomposition Theorem for such maps and for constant coefficients.
Note that if $f$ is not  K\"ahler, see \ci{sa}, then it is not known whether the theorem
holds or not. The counterexamples to the Decomposition Theorem we are aware of 
are not for semismall maps.

\smallskip
We prove Theorem
\ref{fin}: the Decomposition Theorem is equivalent
to the non-dege\-ne\-ra\-cy of a collection of intersection forms that arise
by studying a Mather-Thom-Whitney  stratification of $f$. This result
holds also for more general (i.e. non complex analytic) maps
admitting  stratifications with properties analogous to the complex analytic ones.

\smallskip
Using our Hodge Index Theorem for semismall maps we give a
 new proof of the  Decomposition Theorem for constant coefficients
and for  $X$ and $Y$ projective, Theorem \ref{dtp}. 
This was actually our starting point and 
our main motivation to investigate lef bundles (l.e.f.= {\em Lefschetz
effettivamente funziona}).
A remarkable new feature  we prove
 is that the intersection forms
in questions are not only non-degenerate, but also definite.

\medskip
We believe that the relation between the Decomposition Theorem 
and the intersection forms associated with the strata is illuminating.
Furthermore, the direct proof of the non-degeneracy of these forms, with the 
new additional information on the signatures, neither relying
 on reduction to positive 
characteristic nor on Saito's theory of mixed Hodge modules, sheds 
light on the geometry underlying the decomposition 
theorem and gives some 
indications on its possible extensions beyond the algebraic category.

\bigskip
\noindent
{\bf Acknowledgments.} 
Parts of this work have been done while the first author was visiting
 Hong Kong University in January 2000 and the Max-Planck Institut
f\"ur Mathematik in June 2000  and while the second author was visiting
Harvard University in April 1999 and SUNY at Stony Brook in May 2000.
 We would like to thank
D. Massey, A.J. Sommese and J.A. Wisniewski for useful 
correspondences.

\section{The Hard Lefschetz Theorem for lef line bundles}
\label{hltlb}

In this section we introduce the notion of {\em lef}$\,$ line bundle
on a projective variety. It is a positivity notion 
weaker than  ampleness but stronger than semiampleness
and bigness combined. Lef line bundles satisfy many of the 
cohomological properties  of ampleness. Our main result
is that they also satisfy the Hard Lefschetz Theorem, the Lefschetz Decomposition 
and the Hodge 
Riemann Bilinear relations  on the primitive spaces.
These results are all false, in general,  for line bundles
which are simultaneously generated by their global sections and big.
We also study pencils of sections and obtain more precise statements and a 
second proof of the Hard Lefschetz Theorem. This extra information is 
then used to prove
our  Hodge Index Theorem for semismall maps.

\subsection{Semismall maps and lef line bundles}
\label{smllb}
Let $f:X \to Y$ be a proper holomorphic map. 
For every integer $k$ define $Y^{k}:=
\{ y \in Y \, | \; \dim{ f^{-1} (y)  } = k \}$.
The spaces $Y^{k} $ are locally closed analytic subvarieties
of $Y$ whose disjoint union is $Y$.
If a fiber is reducible, then it is understood that its dimension 
is the highest among the  
dimensions of its components.

\begin{defi}
    \label{defofss}
    {\rm 
We say that a proper holomorphic map $f:X \to Y$
of irreducible varieties is {\em semismall} if 
$\dim{Y^{k}} + 2k \leq \dim{X}$ for every $k$.
Equivalently, $f$ is semismall if and only if
there is no irreducible subvariety $T \subseteq X$ such that
$2 \dim{T} - \dim{ f(T) } > \dim{X}$.
}
\end{defi}

\begin{rmk}
    \label{ssgf}
{\rm A semismall map is necessarily generically finite.}  
     \end{rmk}

From now on we shall assume that semismall maps are proper and 
surjective.

\begin{defi} 
    {\rm 
    We say that a line bundle $M$ on a complex projective variety $X$ is 
    {\em lef} if  a positive multiple of $M$ is generated by its global 
    sections and
 the corresponding morphism onto the image is semismall.
}
\end{defi}

\begin{rmk}
\label{gfip}
{\rm
If the map associated to a multiple of $M$  generated by its global 
sections is  semismall, 
then  the map associated with any other multiple of $M$ generated by 
its global sections is semismall as well. A lef line bundle is nef and big,
 but not conversely.}
\end{rmk}

\begin{pr} 
\label{wl}
(Weak Lefschetz theorem for lef line bundles)
Let $M$ be a lef
line bundle 
on a smooth complex projective variety $X$. Assume that $M$
admits a section   $s \in H^0(X, M)$ whose reduced zero locus
 is a smooth divisor $Y$. Denote by $i: Y \to X$ the inclusion.

\noindent
The restriction map $i^*:H^r(X)\lorw H^r(Y)$ is 
an isomorphism for $r<dimX-1$
and it is
injective for $r=\dim{X}-1$.
\end{pr}
{\em Proof.} The proof can be obtained by a use of the
Leray spectral sequence coupled with the theorem on the cohomological 
dimension of  constructible sheaves on affine varieties.
See, for example, \ci{mig}. See also \ci{e-v}, Lemma 1.2 and Proposition \ref{wlgm}.  
\blacksquare

\medskip
We shall employ the following strenghtening of the Weak Lefschetz Theorem
due to Goresky-MacPherson, \ci{g-m}, II.1.1. In fact their statement
is even stronger.

\begin{pr}
\label{wlgm}
Let $f: X \to Y$ be a semismall holomorphic map from a complex
connected manifold onto a quasi-projective variety $Y$,
and $H \subseteq Y$ be a general  hyperplane section.

\noindent
The restriction map $i^*:H^r(X)\lorw H^r(f^{-1} (H))$ is 
an isomorphism for $r<dimX-1$
and it is
injective for $r=\dim{X}-1$.
\end{pr}

\begin{pr}
\label{ber}
(Bertini theorem for lef line bundles) 
Let $M$ be a lef line bundle 
on a nonsingular  complex projective variety $X$. Assume that $M$
is generated by its global sections. 
Let $W'\subseteq |M| $ be the set of  divisors $Y$  in the linear system of $M$
such that 
$Y$
is smooth and  $M_{|Y}$ is lef.  
Then the  set $W'$ contains a  nonempty and
Zariski open subset $W \subseteq |M|$. 
\end{pr}
{\em Proof.}
By virtue of the standard Bertini Theorem, a generic divisor $D \in 
|M|$
is nonsingular. Let $f: X \to Y$ be the semismall
map associated with $|M|$ and $Y^{k}$  be the locally closed 
subvarieties mentioned above.  The set of divisors containing at least one
among the closed subvarieties $f^{-1}(\overline{Y^{k}})$ is a finite union
of linear proper subspaces of $|M|$. The conclusion follows.
\blacksquare

\begin{rmk}
    \label{kan}
    {\rm  
    The Kodaira-Akizuki-Nakano Vanishing Theorem holds for lef line 
    bundles.  See \ci{e-v}, Theorem 2.4. See also \ci{mig}.}
    \end{rmk}

\begin{rmk}
\label{grariem}
{\rm 
The Kodaira-Akizuki-Nakano Vanishing Theorem for lef line bundles 
coupled with Serre Vanishing Theorem implies that
if $f: X \to Y$ is a semismall holomorphic map from a complex 
projective manifold of dimension $n$ onto a projective variety 
then $R^qf_* \Omega^p_X=0$
for every $p+q >n$. This fact generalizes Grauert-Riemannschneider
Vanishing Theorem and is false for semiample big line bundles.
}
\end{rmk}

\begin{rmk}
    \label{vectbdle}
    {\rm  
    One can define lef vector bundles in terms of the lefness
    of the tautological line bundle and prove that many of  the standard
    properties for ample vector bundles hold  for lef ones.
    }
    \end{rmk}

\subsection{The property HL}
\label{tphl}
Let $X$ be a smooth, compact,  oriented 
manifold of even real  dimension $2n$.
We will use the notation $H^r(X)$ for $H^r(X,\rat)$.
The bilinear form on  $H^{*}(X):= \oplus H^r (X)$
defined by 
$(\alpha,\beta)=\int_{X}{ \alpha\wedge \beta }$  is 
non-degenerate by Poincar\'e
duality.

\noindent
Let  $\omega \in H^2 (X)$. We define   a bilinear form
on $H^{n-r} (X)$ by setting 
$$
 \Psi (  \alpha,\beta ) =
(-1)^{ \frac{(n-r)(n-r-1)}{2}       }
\int_{X}{ \omega^{r}\wedge \alpha 
\wedge \beta },
$$
for every  $0 \leq  r \leq n$.  The form $\Psi$ 
is non-degenerate precisely
when the linear map $L^{r}= L^{r}_{\omega} : H^{n-r} (X) \lorw 
H^{n+r} (X)$, sending
 $\alpha $ to $\omega^{r}\wedge \alpha$,  is an isomorphism.

\bigskip

\begin{defi}
\label{hlr}
{\rm 
We say that $(X,\omega)$
{\em has property $HL_r$} if 
the map $L_{\omega}^r:H^{n-r} (X) \lorw H^{n+r} (X)$ given by the cup product
with $\omega^r$ is an isomorphism. 

\noindent
We say that $(X,\omega)$ has property
$HL$ if it has property $HL_r$ for every $0 \leq r \leq n$.
}
\end{defi}

Note that property $HL_0$ is automatic and that property $HL_n$
is equivalent to $\int_X{\omega^n} \neq 0$.

\noindent
 Define  $ H^{n-r} (X)  \supseteq P^{n-r}= P^{n-r}_{\omega}: =Ker \, 
 L_{\omega}^{r+1}$
and call its elements {\em primitive} (with respect to $\omega$).
The following Lefschetz-type decomposition is immediate.

\begin{pr} 
\label{prim}
Assume that $(X,\omega)$ has property HL.
For every $0 \leq r \leq n$, we have the following ``primitive" decomposition
$$
 H^{n-r} (X) =P^{n-r} \oplus L_{\omega}(H^{n-r-2} (X)).
$$
There is a direct sum decomposition
$$
H^{n-r} (X) = \oplus L_{\omega}^i P^{n-r-2i}.
$$
The subspaces $ L_{\omega}^i P^{n-r-2i}$ are pairwise orthogonal in 
$H^{n-r} (X)$.
\end{pr}

\begin{rmk}
    {\rm 
The
 projection 
of $H^{n-r} (X) $ onto $P^{n-r}$ is given by
$\alpha \to \alpha 
-L_{\omega}(L_{\omega}^{r+2})^{-1}L_{\omega}^{r+1}\alpha$, where
$(L_{\omega}^{r+2})^{-1}$
denotes   the inverse 
to $(L_{\omega}^{r+2}): H^{n-r-2} (X) \to H^{n+r+2} (X)$. 
}
\end{rmk}

\begin{defi}
\label{rhst}
{\rm 
A {\em rational Hodge structure of pure weight $r$}
is  a rational vector space $H$
with a bigraduation of $H_{\comp}=H\otimes {\comp}=\oplus H^{p,q}$  
for $p+q=r$
such that $H^{p,q} =\overline{ H^{q,p} }$.
}
\end{defi}

\begin{defi}
\label{defpol}
{\rm 
A {\em polarization} of the weight  $r$ Hodge structure 
$H$ is a bilinear form $\Psi$ on $H$, symmetric for $r$ even, 
anti-symmetric for $r$ odd, such that its $\comp$-bilinear extension,
to $H_{\comp}$, still denoted by $\Psi$, satisfies :

\noindent
a) the spaces  $H^{p,q}$ and   $H^{s,t}$ are  $\Psi$-orthogonal
whenever either $p \neq t$, or $q \neq s$;

\noindent
b) $i^{p-q} \Psi (\alpha, \overline{\alpha} ) > 0$, for every non-zero
$\alpha \in H^{p,q}$.
}
\end{defi}

\bigskip
Let  $X$ be a nonsingular complex projective variety. For any
ample line bundle $M$ define  
$L_M := L_{ c_{1} (M) }: H^{n-r} (X)  \lorw H^{n-r+2} (X)$. 
Classical Hodge theory gives that 
$\Psi_{M} (\alpha,\beta)=
(-1)^{ \frac{(n-r)(n-r-1)}{2}  }
\int_{X}{ L_M^r(\alpha)\wedge \beta }$
is a polarization of the weight $(n-r)$ Hodge structure  
$P^{n-r}_M=Ker L^{r+1}_M\subseteq H^{n-r} (X)$.

\medskip
We say that $(X, M)$ has property $HL_{r}$ ($HL$, resp.)
if $(X, c_{1} ( M ) )$ has property $HL_{r}$ ($HL$, resp.).

\bigskip
The $HL$ property for a pair $(X, M)$ with $X$ projective and $M$ nef implies 
that $\Psi_{M}$ is a polarization. In fact such a line bundle can be 
written as a limit 
of rational K\"ahler classes and the following proposition applies.

\begin{pr} 
\label{hlpol}
Let $X$ be a compact  connected  complex K\"ahler  manifold  of dimension $n$ and $M$
be a 
 line bundle such that $(X,M)$ has property HL and $c_{1}(M)=\lim_{i 
 \to \infty}{\omega_{i}}$, $\omega_{i}$ K\"ahler.  
 The bilinear form
$\Psi_{M} (\alpha,\beta)=
(-1)^{ \frac{(n-r)(n-r-1)}{2}       }
\int_{X}{ L_M^r(\alpha)\wedge \beta }$
 is a polarization of the weight $(n-r)$
Hodge structure
$P^{n-r}_M=Ker L_M^{r+1}\subseteq H^{n-r} (X)$, for every
$0 \leq r \leq n$.
\end{pr}
{\em Proof.} Since $M$ is nef we can find rational 
K\"ahler classes $\omega_i\in H^2(X)$
such that 
$\lim_{i \to \infty}{\omega_i} = c_1(M)$. 

\noindent
The only thing that needs to be proved is the statement
$i^{p-q}\Psi_M(\alpha, \overline\alpha )>0$ for every non-zero $\alpha\in 
P^{n-r}_M \cap H^{p,q} (X)$.

\noindent
Since the classes $\omega_{i}$ are K\"ahler,
we have the decomposition
$H^{n-r} (X) =P_{\omega_i}^{n-r} \oplus L_{\omega_i}(H^{n-r-2} (X) )$.
Let $\pi_i$ denote the projection onto $P^{n-r}_{\omega_i}$,
$$
\pi_i(\alpha)=\alpha - L_{ \omega_{i} } 
( L^{r+2}_{ \omega_{i} } )^{-1}L_{\omega_{i}}^{r+1}\alpha.
$$
Since $M$ satisfies the HL condition, 
the map  $L_M^{r+2}:H^{n-r-2} (X) \to H^{n+r+2} (X)$
is invertible so that $\lim_{i \to \infty}{ (L_{\omega_{i}}^{r+2})^{-1}} 
= (L_M^{r+2})^{-1}$. Identical considerations hold for the 
$(p,q)$-parts of these invertible maps.
 It follows that, 
if $\alpha\in 
P^{n-r}_M \cap H^{p,q} (X)$, then
$\lim_{i \to \infty}{\pi_i(\alpha)}= \alpha$. 
Since the operators $\pi_{i}$ are of type $(0,0)$, $\pi_{i} (\alpha) 
\in P^{n-r}_{\omega_{i}} \cap H^{p,q} (X)$. Therefore, 
$i^{p-q}\Psi_{\omega_{i}}(\pi_{i}(\alpha), \overline{\pi_{i} (\alpha )} )>0$.
It follows that 
$i^{p-q}\Psi_M(\alpha, \overline\alpha ) \geq 0$. The HL property for 
$M$ implies that  $\Psi_{M}$ is non-degenerate, therefore 
$\Psi_{M}$ is a polarization of $P^{n-r}_{M}$.
\blacksquare

\medskip
The following elementary fact highlights the connection between
the $HL$ property and lef line bundles.

\begin{pr}
\label{nohlnoss}
Let $f:X \to Y$ be a  surjective projective morphism from a nonsingular
projective variety $X$, $A$ be a line bundle on $Y$ and $M:=
f^{*}A$.

\noindent
If $M$ has property $HL$, then $f$
is semismall.
\end{pr}
{\em Proof.} If $f$ is not semismall, then there exists an irreducible
subvariety $T \subseteq X$ such that 
$2 \dim{T} - n > \dim{f(T)}$.
Let $[T] \in H^{2(n- \dim{T})}(X)=H^{n- (2 \dim{T} -n)}(X)$ be the 
fundamental class of $T$. The class $c_{1}(M)^{2 \dim{T} -n}$
can be represented by a $\rat$-algebraic cycle that does not intersect
$T$. It follows that $c_{1}(M)^{( 2 \dim{T} -n )} \cdot [T]=0$, i.e.
$M$ does not satisfy $HL_{2 \dim{T} -n}$.
\blacksquare

\begin{rmk}
    \label{wis}
    {\rm A related fact, dealing with extremal ray contractions, can be found
in 
 \ci{wisn}.  }
     \end{rmk}

\subsection{The Hard Lefschetz Theorem and the signature of intersection
forms}
\label{hlsif}
Our goal is to prove the following extension
of the classical Hard Lefschetz Theorem which  also constitutes a converse
to Proposition \ref{nohlnoss}. At the same time we prove that the 
Hodge-Riemann Bilinear Relations hold on the corresponding primitive 
spaces.
 
\begin{tm}
\label{mainthm}
Let $X$ be a nonsingular complex projective 
variety and  $M$ be a  lef line bundle on $X$.

\noindent
The pair   
$(X,M)$ has  property $HL$. In addition, $\Psi_M$ is a polarization
of $P^{n-r}_{M}= Ker \, L_{M}^{r+1}$.
\end{tm}

\begin{rmk}
    \label{itstopol}
    {\rm 
    Proposition \ref{prim} implies the decomposition of
    the singular cohomology of $X$ into   
    subspaces which are primitive with respect to $M$. It is immediate to check that
    $\dim_{\comp}{P^{l}_{M}}= b_{l} - b_{l-2}$ and $\dim_{\comp}{P^{p+q}_{M} 
    \cap H^{p,q}(X)}= h^{p,q}(X)- h^{p-1, q-1} (X).$
    }
    \end{rmk}

\bigskip
The proof of Theorem \ref{mainthm} can be found in $\S$\ref{firstproof}.
A second proof, based on 
 a variation of Deligne's proof of the Hard Lefschetz
Theorem, see \ci{weil2}, can be found in $\S$\ref{proofhl}.  

\subsubsection{Proof of Theorem \ref{mainthm}}
\label{firstproof}
We use induction on the dimension of $X$.
Note that the case $n=1$ is classical, for $M$ is then 
necessarily ample. The statement
is invariant under taking  non-zero positive powers of the line bundle $M$
so that we shall assume from now on that $M$ is generated by its 
global sections.

\medskip
In this section $X$ denotes a nonsingular 
complex  projective variety
of dimension $n\geq 2$, and $M$
a lef line bundle on $X$ generated by its global sections.
Let $s \in H^0(X, M)$ be a section with smooth 
zero locus 
 $Y$ such that $M_{|Y}$ is lef. Such a section exists by 
 Proposition \ref{ber}. Note that  $Y$ is necessarily connected
 by Bertini Theorem. 
 Denote by
$i: Y \to X$ the inclusion.

Let $\hat{L}_M^r: = (L_{ M_{| Y} })^{r}: H^{n-1-r} (Y) \to 
H^{n-1+r } (Y)$.
We record, for future use,  the following elementary fact that
 follows at once  from  the projection formula
and Poincar\'e Duality.

\begin{lm}
\label{gy} 
$L_M^r=i_* \circ \hat{L}_M^{r-1}\circ i^*$.
\end{lm}

\begin {lm}
\label{yhl1} 
If $(Y,M_{|Y})$ has  property HL, then $(X,M)$ has
properties $HL_r$ for $r=0$ and $ 2 \leq r \leq n$.

\noindent
The pair $(X,M)$  has  property $HL_1$
if and only if 
the restriction of the intersection form on $H^{n-1}(Y)$ to the subspace
$i^*H^{n-1}(X)$ is non degenerate.

\noindent
If $(Y,M_{|Y})$ has  property HL, then $(X,M)$ has
property $HL$ 
if and only if 
the restriction of the intersection form on $H^{n-1}(Y)$ to the subspace
$i^*H^{n-1}(X)$ is non degenerate.
\end{lm}
{\em Proof.} 
The case $r=0$ is  trivially true.
Let  $2 \leq r \leq n $.  
The first assertion follows from Proposition \ref{wl},
the assumption that property $HL$ holds for $(Y,M_{|Y})$ and 
 Lemma \ref{gy}.

\noindent
Let us prove the second statement.
The map $L_M: H^{n-1}(X) \lorw H^{n+1}(X)$ is the composition of the 
injective map $i^*:H^{n-1}(X) \lorw H^{n-1}(Y)$ and its surjective
transpose $i_*:H^{n-1}(Y) \lorw H^{n+1}(X)$. It is an isomorphism iff
$i^*H^{n-1}(X) \cap Ker \,i_* = \{\,0\,\}$, i.e. iff 
$i^*H^{n-1}(X) \cap (i^*H^{n-1}(X))^{\perp}= \{\,0\,\}$. This
latter statement is equivalent
to the non-degeneracy of the restriction of the intersection form.
The third assertion follows. 
\blacksquare

\bigskip
We now prove Theorem \ref{mainthm}.
By Lemma \ref{yhl1} it is enough to show the statement of non-degeneracy
on $i^*H^{n-1}(X)$. Conisder $i_*: H^{n-1}(Y) \to H^{n+1}(X)$.
We have $ (i^* H^{n-1}(X) )^{\bot}= Ker \, i_* \subseteq 
P_{ M_{|Y } }^{n-1}(Y)$. By induction this last space is polarized
by the intersection form. In particular, the intersection form in non-degenerate
on $(i^* H^{n-1}(X) )^{\bot}= Ker \, i_*$ so that it is non-degenerate
on $i^* H^{n-1}(X)$.
It follows  that $(X,M)$ has  property HL. 
We conclude by Proposition \ref{hlpol}.
\blacksquare

\subsubsection{The use of pencils}
\label{pencils}
In this section we place ourselves in the same situation 
as in $\S$\ref{firstproof}. The goal is to prove Theorem
\ref{semis}.

\bigskip
Choose a  pencil $l \subseteq |M|$ 
which meets the set $W$ of  Proposition
\ref{ber}. There are two elements
$u_0$ and $u_{\infty} \in  l \cap W$ such that the corresponding 
nonsingular divisors
meet transversally along a nonsingular subvariety $A$. 
By blowing up $X$ along $A$ we obtain a nonsingular variety $X'$ together
with a proper and flat  morphism $f: X' \to l \simeq \pn{1}$. 
Denote by $U$ the Zariski open and dense subset $l \cap W$. The morphism
$f$ is smooth over $U$.
For every $u \in U$, the fiber $f^{-1}(u)$ is naturally isomorphic
to the divisor $Y_u$ corresponding to $u \in U \subseteq l$ and, 
in particular, $M_{|Y_u}$
is lef for every $u \in U$.
We denote $Y_{u_0}$ simply by $Y$ and we denote the strict transform
of $Y$ on $X'$ by $Y'$. Clearly, $Y$ and $Y'$ are naturally isomorphic.
Let $ j: U \to \pn{1}$ denote the open imbedding.
The sheaves $j^{*}R^q f_* \rat_{X'}$ 
are local systems on $U$, for every $q \geq 0$,
corresponding to the 
  $\pi_1(U,u_0)$-module
$ H^q(Y')$. By abuse of notation, we drop the base point
$u_0$.

\medskip
The following is well-known, e.g. \ci{looj}, Lemma 5.3:
\begin{lm}
    \label{inv=sect}
    Let
${\cal V}$ be   local system on $U$ associated with a finite 
dimensional
$\pi_1(U)$-mo\-du\-le $V$.

\noindent
We have that 
$H^0(\pn1,j_*{\cal V})=H^0(U,{\cal V})= V^{\pi_1 (U)}$,
the submodule of invariants.
\end{lm}

\bigskip
We now invoke the following well-known  result of Deligne's

\begin{pr}
    \label{iandj}
Let $i':Y' \lorw f^{-1}(U)$ and 
$j':f^{-1}(U)\lorw X'$ be the corresponding closed and open imbeddings. 
For every $q$ we have
$$
 i'^* H^q(f^{-1}(U)) = (j'\cdot i')^* H^q(X') 
\subseteq  H^q(Y').
$$ 
\end{pr}
{\em Proof.} See \ci{ho2},  Corollaire 3.2.18.
 \blacksquare

\bigskip

In what follows we do not want to use Theorem
\ref{mainthm}, so that our study will also yield a second proof of that result.
\begin{pr}
\label{ler}
Assume that  Theorem \ref{mainthm} holds
in dimension $n-1$.

\noindent
The Leray spectral sequence for the restriction of $f$ to 
$f^{-1}(U)$ degenerates
at $E_2$. In particular, the morphism
$$
H^q(f^{-1}(U)) \lorw H^0(U,R^qf_* \rat_{X'})=H^q(Y)^{\pi_1 (U)}
$$
is surjective for every $q$.
\end{pr}
{\em Proof.} By virtue of the assumptions and of the choice of $U$, 
we can apply the degeneration criterion in   \ci{dess}. 
\blacksquare

\begin{cor}
\label{inv}
Assume that  Theorem \ref{mainthm} holds
in dimension $n-1$.
The image of  $H^q(X')$ in $H^q(Y')$ by the restriction map
is $H^q(Y')^{\pi_1 (U)}$ for every $q$.
\end {cor}
{\em Proof.} The restriction map 
$$
i'^*:H^q(X') \lorw H^q(Y')
$$
factors as follows
$$
H^q(X') \lorw H^q(f^{-1}(U))
\lorw H^0(U,R^q f_*\rat) \lorw H^q(Y')^{\pi_1(U)}
\lorw H^q(Y'). 
$$
The statement follows from Proposition  \ref{iandj}  and Proposition 
\ref{ler}.
\blacksquare

\bigskip
We now  compare 
$H^q( X) \lorw H^q(Y)$ with
$H^q(X') \lorw H^q(Y') $.

\begin{pr}
Let $i':Y'\lorw X'$ and $Bl_A: X' \lorw X$ the blowing up map.
We have that $i'^* \cdot Bl_A^*$ is injective. In addition, if $\, q\leq n-1$,
then  
${i'}^*H^q(X')={i'}^* \cdot Bl^*_A H^q(X) 
\subseteq H^q(Y')$.   
\end{pr}
{\em Proof.} See \'Expos\'e XVIII Corollaire 5.1.6 in SGA VII,2.
\blacksquare

\bigskip

Again, we do not want to use Theorem \ref{mainthm}.
Recall that $Y$ and $Y'$ are naturally identified.
\begin{tm}
\label{semis}
Assume that  Theorem \ref{mainthm} holds
in dimension $n-1$.
For $q\leq \dim{X} -1$ 
the map    $i'^* \cdot Bl^*_A: H^q(X) \lorw H^q(Y')$
is injective and its image is precisely the subspace $H^q(Y')^{\pi_1 (U)}$
of monodromy invariants.
In addition,
the $\pi_1(U)$-representation on $H^{n-1}(Y')$ is semisimple.
\end{tm}
{\em Proof.} We only  need to prove the last statement. 
The restriction of $Bl_A^*M$ to $f^{-1}(u)$ is lef
for every $u \in U$.  By virtue of the assumptions and 
of Theorem \ref{hlpol},
this defines a polarization on all the primitive subspaces 
$P^q_M(u) \subseteq H^q(f^{-1}(u))$. It follows   that 
$H^{n-1}(f^{-1}(u))=\oplus L_{M}^i P^{n-1-2i}_M(u)$.
The local systems
$P^{n-1-2i}_M(u)$ define a variation of polarized Hodge structures. 
We conclude by the semisimplicity theorem of Deligne 
\ci{ho2},  Theorem  4.2.6 and footnote 1.
 \blacksquare

\subsubsection{A second proof of  Theorem \ref{mainthm}}
\label{proofhl}
We prove the theorem by induction on $\dim{X}=n$. The statement is true 
when
$\dim{X} =1$, for then $M$ is necessarily ample.

\noindent
Let us assume that the statement is true in dimension $n-1$.
We replace, without loss of generality, $M$ by a positive power which is
generated by its global sections. 

\noindent
We  fix a pencil $l$ as in $\S$\ref{pencils}.
By virtue of Theorem \ref{semis},  we can find a $\pi_1 (U)$-invariant complement $V$
to $H^{n-1}(Y')^{\pi_1 (U) }$ in $H^{n-1}(Y')$. The Poincar\'e duality pairing
on $H^{n-1}(Y')$ is $\pi_1 (U)$-invariant, therefore it restricts to 
a non-degenerate pairing on  
 $H^{n-1}(Y')^{\pi_1  (U)} $.


 \noindent
We conclude that $(X,M)$ has the property HL by virtue of 
the third part of Proposition \ref{yhl1}. 
Finally,  we conclude by Proposition \ref{hlpol}.
\blacksquare

\subsection{The Hodge Index Theorem for semismall maps}
\label{hith}
Let us record the following consequence of 
Theorem \ref{semis} and of Proposition \ref{wlgm}.
Together with Remark \ref{anynumber}, they are
a higher dimensional analogue of Grauert
contractibility test for curves on surfaces.

\begin{tm}
\label{fai} 
(Hodge Index Theorem for semismall maps)
Let $f: X \to Y$ be a  semismall  map from 
a nonsingular complex projective variety of even dimension 
$n$ onto a projective variety $Y$ and $ y \in Y $ be  a  point such that
 $\dim { f^{-1}(y) } = \frac{n}{2}$. Denote by $Z_{l}$, $1 \leq l \leq r$, 
 the irreducible components
 of maximal dimension of $f^{-1} (y)$.
 
 \noindent
 Then the cohomology classes $[Z_l ] \in H^n(X)$ are linearly independent
 and the symmetric matrix $  (-1)^{ \frac{n}{2} } || Z_l \cdot Z_m || $
 is positive definite.
\end{tm}
{\em Proof.}
 Choose  a general pencil  $l$ of divisors on $X$   as in $\S$\ref{pencils}
 with the additional property that its general members
 satisfy the conclusion of Proposition \ref{wlgm} for 
 $X \setminus f^{-1} (y)$. Let $V \subseteq l$ be the Zariski dense 
 open subset such that the members  $X_{v}$ are nonsingular, 
 $X_{v}\cap f^{-1} (v)  =\emptyset$ and Proposition \ref{wlgm}
 holds for $X \setminus f^{-1} (y)$ and $X_{v}$.
 We have that the restriction map $H^{n-1} (X \setminus f^{-1} (y) )
 \to H^{n-1} (X_{v} )$ is injective for $v \in V$.
 By virtue of  Theorem \ref{semis}, we have the isomorphism
 $H^{n-1}(X) \to H^{n-1}(X_v)^{\pi_1 (V) }$, for every $v \in V$.
 This isomorphism factors as follows
 $$
 H^{n-1} (X) \to H^{n-1} (X \setminus f^{-1}(y)) \stackrel{a}\lorw
  H^{n-1}(X_v)^{\pi_1 (V)} \subseteq  H^{n-1}(X_v).
  $$
Hence the injective  map $a$ is 
surjective. 

\noindent
It follows that the natural cycle class map $H^{BM}_n (f^{-1} (y) )  \simeq 
H^n (X, X \setminus f^{-1} (y) ) \to H^n (X)$ is injective.
Its image is contained in the primitive space associated with any line bundle
on $X$ pull-back of an ample line bundle on $Y$.
The statement follows from  Theorem \ref{mainthm}.
\blacksquare

\begin{cor}
    \label{purity}
Let $f: X \to Y$ and $y$ be as above.
    The mixed Hodge Structure on $H^{n-1} (X \setminus f^{-1} (y) )$ is 
pure of weight $n-1$.    
    \end{cor}
{\em Proof.} It follows from the surjectivity of the natural map
$H^{n-1} (X) \to H^{n-1} (X \setminus f^{-1} (y) )$.
\blacksquare

\begin{cor}
    \label{issing}
Let  $f : X \to Y$ be a birational semismall map from
a nonsingular   quasi projective complex variety of even dimension $n$ onto a 
quasi projective complex variety with an  isolated singularity $y\in Y$ 
such that $f$ is an isomorphism outside $y$ and $\dim{f^{-1} (y)}= 
\frac{n}{2}$. Then the conclusions of Theorem \ref{fai} hold. 
    \end{cor}
{\em Proof.} 
One finds a semismall projective completion $f'': X'' \to Y''$ of $f$
to which we  apply Theorem \ref{fai}. Since the bilinear form on the 
fibers in non-degenerate, the cycle classes of the fundamental classes of the 
fibers stay 
independent in $H^{n} (X)$.
\blacksquare

\begin{rmk}
\label{anynumber}
{\rm 
The same proof shows that the statements of Theorem \ref{fai}
and Corollary \ref{issing} holds unchanged if we consider
any finite number of fibers over points as $y$ above.
}
\end{rmk}

\begin{rmk}
     \label{grauerttype}
     {\rm 
 If  $(Y, y) $ is a germ of a normal  complex space of dimension two, 
and $f:X \to Y$ is a resolution of singularities, then 
 Grauert Contractibility Criterion, see \ci{laufer}
Theorem 4.4,  implies that the form in question 
 is negative definite and, in particular, it is non-degenerate.
 }
 \end{rmk}

The following is a natural question. A positive answer would 
yield a proof of the Decomposition Theorem for semismall holomorphic
 maps from  complex manifolds and for constant coefficients; see
Theorem \ref{fin}.

\begin{??}
    \label{signs}
{\rm 
Let $f: V \to W$ be a proper holomorphic semismall map from a complex 
manifold of  even dimension $n$ onto an  analytic space
$Y$. Assume that the fiber $f^{-1}(w)$ 
over a point  $w \in W$ has dimension $\frac{n}{2}$. Is the 
intersection form on $H^{BM}_{n}(f^{-1}(w))$ non-degenerate? Is it  
$(-1)^{\frac{n}{2}}$-positive definite?
}
     \end{??}

\section{The topology of semismall maps}
\label{tssm}
We now proceed to a study of holomorphic semismall maps from 
a complex manifold. First we need to prove  Proposition \ref{splits}, a simple  splitting criterion in derived categories for which we could not find a reference. We study the topology of these maps by attaching one stratum at the time. In doing so a symmetric bilinear form
emerges naturally; see Proposition \ref{rocom} and Lemma \ref{itssymm}.
We then prove that the Decomposition Theorem for these maps and for constant coefficients is equivalent to the non-degeneration of these forms; Theorem
\ref{fin}. Finally, we give a proof of the Decomposition Theorem
when the domain and target are projective, Theorem \ref{dtp}.
 A new feature that we discover is that the forms are definite
by virtue of our Hodge Index Theorem for semismall maps.

\subsection{Homological algebra}
\label{homalg}
 Let ${\cal A}$ be an abelian category with enough injectives, 
 e.g. sheaves of abelian groups on a topological space,
 and $C({\cal A} )$ be the associated category of complexes.
 Complexes and morphisms can be truncated.
  Given an integer $t$, we have two types
 of truncations: $\tau_{\leq t} A$ and $\tau_{\geq t} A$. The former is defined
 as follows $(\tau_{\leq t} A)^i:=A^i$ for $i\leq t-1$, $(\tau_{\leq t} 
 A)^t:=
 Ker (A^t \to A^{t+1} )$, $( {\tau_{\leq t} A} )^i:= \{0 \}$ for $i >t$.
 The latter is defined
 as follows $(\tau_{\geq t} A)^i:=\{0 \} $ for $i\leq t-1$, $( \tau_{\geq 
 t} A )^t:=
 Coker (A^{t-1} \to A^{t} )$, $( \tau_{\geq t} A )^i:= A^i$ for $i >t$.
 Let $h:A \to B$ be a morphism of complexes.
 The truncations $\tau_{\leq t} (h): \tau_{\leq t} A \to \tau_{\leq t} B$
 and  $\tau_{\geq t} (h): \tau_{\geq t} A \to \tau_{\geq t} B$ are defined
 in the natural way. 
 The operations of truncating complexes and morphisms of complexes induce 
 functors in the derived category $D({\cal A} )$.
 
 If $A$ is a  complex acyclic  in degrees $l\neq t$ for some integer 
 $t$, i.e. if  $\tau_{\leq t} A \simeq \tau_{\geq t} A$,  then  $A \simeq
 {\cal H}^{t} (A) [-t]$.
 
 The cone construction for a  morphism of complexes  $h:A \to B$ gives 
 rise, in a non-unique way, 
 to a diagram of morphism of complexes
 $ A\stackrel{h}\to B \to M(h) \stackrel{[1]}\to A[1]$.   
 A  diagram of 
 morphisms $X \to Y 
 \to Z \stackrel{[1]}\to X[1]$ in $D({\cal A})$ 
 is called a distinguished triangle if it 
 is isomorphic to a diagram arising from a cone. 
 
 A morphism $h: A \to B$ in $D({\cal A})$
 gives rise to a distinguished triangle 
$ A  \stackrel{h}\to B \to C \to A[1]$.
 If $h=0$, then $C \simeq A[1] \oplus B$ and  the induced 
 morphism $A[1] \to A[1]$ is an isomorphism.

 A morphism $h: A \to B$ in the derived category gives a collection of morphisms in  cohomology
${\cal H}^l(h): {\cal H}^l(A) \to {\cal H}^l (B)$. A distinguished triangle 
$A \to B \to C \stackrel{[1]}\to A[1]$ gives rise to a cohomology long exact sequence:
$$
\ldots \;  {\cal H}^l(A) \to {\cal H}^l(B) \to {\cal H}^l(C) \to {\cal 
H}^{l+1}(A) \;
\ldots
$$

A non-zero morphism  $h: A \to B$ in the derived category 
may nonetheless induce the zero morphisms between all cohomology groups. 
However, we have the following simple and standard
 \begin{lm}
\label{hom}
Let $t$ be an integer and $A$ and $B$ be two complexes such 
that $A \simeq \tau_{\leq t}A$ and $B\simeq \tau_{\geq t} B$. 
Then the natural map
$Hom_{ D ({\cal A}) } (A,B)  \to 
Hom_{\cal A} ( {\cal H}^t (A)  , {\cal H}^t ( B)  )$ 
is an isomorphism of abelian groups.
\end{lm}
{\em Proof.}  It is enough to replace $B$ by an injective resolution
placed in  degrees no less than $t$. 
\blacksquare
 
\bigskip
 We shall need the following elementary splitting criterion.
 \begin{pr}
 \label{splits}
 Let $ C \stackrel{u}\lorw A \stackrel{v}\lorw B \stackrel{[1]}\lorw C[1]$ be a 
 distinguished triangle and $t$
 be an integer
 such that $A \simeq  \tau_{ \leq t} A$ and $C \simeq  \tau_{ \geq t} C$.
 
 \noindent
 Then
 ${\cal H}^t (u) : {\cal H}^t(C) \to {\cal H}^t (A)$ is an isomorphism
 iff
 $$
 A \simeq \tau_{\leq t-1}B \,  \oplus \,  {\cal H}^t(A) [-t]
 $$
 and the map $v$ is the direct sum of the natural map
 $\tau_{\leq  t-1} B \to B$ and the zero map.
 \end{pr}
{\em Proof.} Assume that ${\cal H}^t (u)$ is an isomorphism.
Apply  the functor $Hom (A, -)$
to the distinguished triangle 
$\tau_{ \leq t-1 } B  \stackrel{ \nu_{ t-1}  } \lorw   
  \tau_{ \leq t} B  \stackrel{\pi}\lorw {\cal H}^t (B) [-t] \stackrel{[1]}\lorw
   \tau_{ \leq t-1} B  \,[1]$
and we get the following exact 
sequence:
$$
\ldots \to  Hom^{-1}  ( \tau_{ \leq t} A , {\cal H}^{t} (B) [-t] )
\to Hom^{0}  ( \tau_{ \leq t} A , \tau_{ \leq t-1} B ) 
\to   
$$
$$
Hom^{0}  ( \tau_{ \leq t} A , \tau_{ \leq } B ) \to
 Hom^{0}  ( \tau_{ \leq t} A ,  {\cal H}^{t} (B) [-t] ) )  \to 
 \ldots
$$
Since ${\cal H}^{t} (B) [-t] $ is concentrated in 
degree $t$, $Hom^{-1}  ( \tau_{ \leq t} A , {\cal H}^{t} (B) [-t] ) 
=\{0\}$.  The morphism ${\cal H}^{t} (v) = 0$, for ${\cal H}^{t} (u)$ 
is  surjective.

\noindent
It follows that there exist a unique lifting $v'$ of $\tau_{\leq t} 
(v)$, i.e. there exists a unique $v': A \to \tau_{\leq t-1} B$ such that
$\tau_{\leq t} (v)= \nu_{t-1} \circ v'$.

\noindent
We complete $v'$ to a distinguished  triangle:
$$
\tau_{ \leq t} A \stackrel{v'}\lorw \tau_{ \leq t-1} B \stackrel{v''}\lorw
M(v') \stackrel{[1]}\lorw \tau_{ \leq t} A  \, [1] \;.
$$
By degree considerations, the morphism ${\cal H}^l (v')=0 $ for $l \geq t$.
Since $v'$ is a lifting of $\tau_{ \leq t  } (v)$, the morphism ${\cal H}^l (v')$
is an isomorphism for $l\leq t-1$ and it is the zero map for $l \geq t$.
This implies that $M(v' ) \simeq  {\cal H}^t (A) [-t +1]$ and that ${\cal H}^{t-1} (v'')=0$.
By virtue of Lemma \ref{hom}, we get that $v'' =0$.

\noindent
The desired splitting follows. The converse can be read off the long 
exact cohomology sequence.
\blacksquare

\subsection{The bilinear forms associated with relevant strata}
\label{bfars}
\label{bfa}
Let $f: X \to Y$ be a proper holomorphic semismall map with $X$ nonsingular
connected 
of dimension $n$. 
Let us summarize the results from 
stratification theory (cf. \ci{g-m}, Ch. 1) that we shall need in the sequel.
They are based essentially on Thom First Isotopy Lemma.

\bigskip
There exists a collection of disjoint locally closed and {\em connected}
analytic 
subvarieties $Y_{i} \subseteq Y$ such that:

\medskip
a) $Y= \coprod_i{Y_i}$ is a Whitney stratification of $Y$.

b) $Y_i \cap \overline{Y_j} \neq \emptyset$ iff $Y_i \subseteq 
\overline{Y_j}$.

c) the induced maps $f_i: f^{-1}(Y_i) \to Y_i$ are stratified submersions; 
in particular they are topologically locally trivial fibrations.




\bigskip
We call such data {\em a stratification of the map $f$}.

 \begin{defi}
     \label{strel}
{\rm A stratum $Y_{i}$ is said to be {\em relevant} if 
$2 \dim{ f^{-1}( Y_{i}) }-  \dim{Y_{i}}=n$. Let $I' \subseteq I$
be the set of indices labeling relevant strata.
}    
 \end{defi}

Let $i \in I$ be any index 
and $d_{i}:= \dim{Y_{i}}$. 
Define
${\cal L}_{i} : = ( R^{ n- d_{i} }  f_* \rat_{X} )_{|Y_{i}}$.
It is  a  local system on $Y_{i}$. 

 \begin{rmk}
     \label{ifnotrelzero}
     {\rm 
 If $Y_{i}$ is not relevant, then
 ${\cal L}_{i}$ is the zero sheaf.
 If $Y_{i}$ is relevant, then the stalks  $( {\cal L}^{*}_{i} )_{ y_{i} }
  \simeq H^{BM}_{n-d_{i}} ( f^{-1} (y_{i} ) )$
 of the dual local system are generated exactly by 
  the fundamental classes of the irreducible and reduced components of maximal 
  dimension of the fiber over $y_{i}$. 
     }
     \end{rmk}
     
 The following  is elementary and holds also when the stratum is not 
relevant when we consider the local system dual to the one
generated by the components
of maximal dimension of the fibers.
 \begin{lm}
 \label{lsss}
 The local system ${\cal L}_i$ splits as a direct sum ${\cal L}_i \simeq
 \bigoplus_{j=1}^{m_i}{ {\cal L}_{ij} }$ of irreducible local sub-systems.
 \end{lm}
 {\em Proof.}
  It is enough to show the   statement for ${\cal L}_{i}^{*}$, $i \in 
 I'.$
   Going around a loop in $Y_{i}$ has the effect
   of permuting the elements of  the basis of Remark 
   \ref{ifnotrelzero}. The associated 
monodromy representation factors 
 through  a finite  symmetric group so that  
  it splits into a sum of irreducibles.
 \blacksquare

 \bigskip
Let $S:=Y_{i}$, $d: =\dim{S}$ and ${\cal L}_S:= {\cal L}_i$.
We now proceed to associating with $S$ a symmetric bilinear form on 
the local system ${\cal L}_{S}^{*}$.

Let $s \in S$ and choose a small-enough euclidean neighborhood $U$ of $s$ in 
$Y$ such that a) $S':= S \cap U$ is contractible and b) the restriction
$ i^{*}: H^{n-d} (f^{-1}(U)) \to H^{n-d} ( f^{-1} (s)  )$ is an 
isomorphism.

Let $F_{1}, \ldots, F_{r}$ be the irreducible and reduced components of maximal 
dimension of $f^{-1} (S')$. By virtue of a) above and of the 
topological triviality over $S'$, the intersections 
$f_{j} : = f^{-1} (s) \cap F_{j}$ are exactly the irreducible and reduced
components of maximal dimension of $f^{-1} (s)$. The analogous 
statement  is true for 
every point $s' \in S'$ and the components for the point
$s$ can be canonically identified with the ones of $s'$.
The specialization morphism
$i^{!}_{s}: H^{BM}_{n+d} (f^{-1} (S')) \to H^{BM}_{n-d} (f^{-1} (s))$,  
associated with the regular imbedding $i_{s}:  \{s \} \to S'$, sends
the fundamental class of a component   $F_{l}$ to the fundamental class
of the corresponding $f_{l}$ and it is an isomorphism; see
\ci{fulton}, Ch. 10.
We have the following sequence of maps:
$$
H^{BM}_{n-d} (f^{-1} (s) ) \stackrel{(i^{!}_{s})^{-1}}\lorw
H^{BM}_{n+d} (f^{-1} (S')) \stackrel{ (\cap \mu_{f^{-1} (U)} )^{-1}   } 
\lorw
H^{n-d} (   f^{-1} (U), f^{-1} (U \setminus S')  ) 
$$
$$
\stackrel{nat}\to
H^{n-d} ( f^{-1} (U)   ) \stackrel{i^{*}}\lorw H^{n-d} ( f^{-1} (s)  )
\stackrel{\kappa}\lorw H^{BM}_{n-d} (f^{-1} (s) )^{*}.
$$
The second map is the inverse to the isomorphism given by capping 
with the fundamental class $\mu_{f^{-1} (U)}$ (cf. \ci{iv}, IX.4).
The third map is the natural map in relative cohomology. The fourth 
map is an isomorphism by virtue of condition b) above.
The map $\kappa$ is an isomorphism by the compactness of $f^{-1} (s)$.

We denote the composition, which is independent of the choice of $U$:

$$
\rho_{S,s}: H^{BM}_{n-d} (f^{-1} (s) ) \lorw  H^{BM}_{n-d} (f^{-1} 
(s) )^{*}.
$$
We have that
$$
\rho_{S,s} (f_{h} ) (f_{k} ) = \deg{F_{h}\cdot f_{k}},
$$
where the refined intersection product takes place in $f^{-1} (U)$ 
and has values in $H^{BM}_{0} (f^{-1} (s))$.

Since the map $f$ is locally topologically trivial along $S$, the 
maps $\rho_{S,s}$ define a map of local systems
$$
\rho_{S} = \rho: {\cal L}_{S}^{*} \lorw {\cal L}_{S}.
$$

\medskip
We record the following fact for future use.
\begin{pr}
\label{rocom}
If $S$ is not relevant, then $\rho_S$ is the zero map between 
trivial local systems.
Let $s \in S$ be a point. The map $\rho_{S,s}$ is an isomorphism
iff the natural map $r_k: H^k(f^{-1}(U)) \to H^k(f^{-1}(U\setminus S'))$
is an isomorphism for every $k \leq n-d-1$, iff  
the natural map $s_k: H^k(f^{-1}(U\setminus S')) 
\to H^{k+1}(f^{-1}(U),f^{-1}(U \setminus S'))$
is an isomorphism for every $k \geq n-d$.
   \end{pr}
{\em Proof.} The domain and the range of $\rho_{S,s}$ are dual to each other.
The statement follows from the relative cohomology sequence for the pair
$(f^{-1}(U),f^{-1}(U \setminus S'))$, the isomorphisms
$H^k(f^{-1}(U)) \simeq H^k(f^{-1}(s))$, 
$H^k(f^{-1}(U),f^{-1}(U \setminus S')) \simeq H_{2n-k}^{BM}(f^{-1}(S'))$
and the fact that $\dim{f^{-1}(s)}  \leq \frac{n-d}{2}$, 
$\dim{f^{-1}(S)}\leq \frac{n+d}{2}$.
\blacksquare

\bigskip
Since $f_i: f^{-1}(S) \to S$ is a stratified submersion, 
given any point $s\in S$, we can choose
an analytic normal slice $N(s)$  to $S$ at $s$ such that
$f^{-1} ( N(s) )$ is  a locally closed complex  submanifold of $X$
of dimension $n-d$.
We now use this fact to express the map $\rho_{S,s}$ in terms of
the refined intersection pairing on $f^{-1} (N(s))$.

\begin{lm}
    \label{itssymm}
 If $s \in S$, then $\rho_{S,s} (f_{h} ) (f_{k} ) = \deg{f_{h} \cdot 
 f_{k}}$, where the refined intersection product on the r.h.s. takes 
 place
 in $f^{-1} (N (s) )$ and has values in $H^{BM}_{0} ( f^{-1} (s) )$.
In particular, the map $\rho_{S}: {\cal L}_{S}^{*} \lorw {\cal L}_{S}$ is symmetric.
\end{lm}
 {\em Proof.} 
Since $f_i: f^{-1}(S) \to S$ is a stratified submersion
and $N(s)$ is a normal slice to $S$ at $s$, $F_j$ meets
$f^{-1}(N(s))$ transversally at the general point of $f_j$.
It follows that the refined intersection 
 product $f^{-1} (N (s) ) \cdot F_{j}$ is the fundamental class
 of $f_{j}$ in $H^{BM}_{n-d} (f^{-1} (s) )$. The result follows by 
 applying \ci{fulton}, 8.1.1.a) to the maps
 $f^{-1} (s) \to f^{-1} (N (s) ) \to f^{-1} ( U )$.
 \blacksquare

\subsection{Inductive study of  semismall analytic maps}
\label{issam}
\bigskip

Let $f:X \to Y$ and $\{ Y_j \}$, $i\in I$, $S$ be as in section \ref{bfa}.
 We  assume, for simplicity, the 
 $Y_{i}$ to be connected and $I$ to be finite. 
There is no loss of generality, for strata of the same 
 dimension do not interfere with each other from the point of view of 
 the  analysis that follows  and could be treated simultaneously.  
 As usual, we define  a partial order
on the index set $I$ by setting $i\prec j$ iff $Y_{i}\subseteq \overline {Y_{j} }$.
We fix a total order 
$I= \{i_{1}< \ldots <i_{\iota}\}$ which is compatible with the 
aforementioned partial order
and define the open  sets $U_{\geq i}:=\coprod _{j\geq i}Y_{j}$.
Similarly, $U_{> i}:=\coprod _{j>i}Y_{j}$.
Let  $\alpha_i: U_{>i} \to U_{\geq i}$ 
be the open imbedding.
We can  define the intermediate extension  
of a complex of 
sheaves $K^{\bullet}$ on $U_{>i}$ to a complex of sheaves 
on 
$U_{\geq i}$ by setting
$$
{\alpha_i}_{!*}K^{\bullet}= \tau_{\leq - dimY_{i}-1 } R{\alpha_{i}}_{*} K^{\bullet}. 
$$
See \ci{bbd}.
The construction is general and can be iterated so that one can form 
the intermediate extension of a complex of sheaves on any $Y_{i}$
to a complex on $\overline{Y_{i}}\cap U_{>j}$ for $j<i$.
In particular, let ${\cal L}$ be a local system on $Y_{i}$. The 
intermediate extension of ${\cal L}[\dim{ Y_{i} }]$ to $\overline 
Y_{i}\cap U_{>j}$ for $j<i$ is called the {\em intersection cohomology 
complex} associated with ${\cal L}$ and is denoted by $ IC_{\overline 
Y_{i}\cap U_{>j}}({\cal L})$.

 \begin{defi}
    \label{defdecthm}
    {\rm 
 Let $f:X \to Y$ be a proper holomorphic  semismall map  from a nonsingular 
connected complex manifold $X$ of dimension $n$. 
We say that {\em the Decomposition Theorem holds for $f$} if there is
an isomorphism
  $$
  Rf_* \rat_{X} [n] \simeq \bigoplus_{ k \in I' }{
  IC_{\overline{Y_{k}}}   ({\cal  L}_{k} )} \simeq 
  \bigoplus_{ k \in I' }\bigoplus_{m=1}^{m_k}{
  IC_{\overline{Y_{k}}}   ({\cal  L}_{km} )},
  $$
  where the ${\cal L}_{km}$ are as in Lemma \ref{lsss}.
}
\end{defi}    
  \begin{rmk}
      \label{whenitholds}
{\rm 
The Decomposition Theorem holds, in the sense defined above, for 
$X$,  $Y$ and $f$ algebraic (cf. \ci{bbd}) and for $f$ a K\"ahler 
morphism (cf. \ci{sa}). In both cases, a far more general 
statement  holds. As observed in \ci{b-m}, $\S$1.7, in the case of semismall
maps these results
can be expressed in the convenient form of Definition \ref{defdecthm}.}
\end{rmk}

\bigskip
 We now proceed to show  that the non-degeneracy of the forms
 $\rho_{S}$ associated with the strata $Y_{i}$
implies the Decomposition Theorem.

\bigskip
Recall that
$I' \subseteq I$ is the subset labeling relevant strata.
For ease of notation set
$$
V : = U_{> i}, \qquad V' : = U_{\geq i},  \qquad  S:= Y_{i}
$$
and let 
$$
V \stackrel{\alpha}\lorw V' \stackrel{\beta}\leftarrow S
$$
be the corresponding open and closed imbeddings.

\begin{tm}
\label{fin}
Assume that the Decomposition Theorem holds over $V$.
The map $\rho_S$  is an isomorphism iff the Decomposition Theorem holds over $V'$ 
and the corresponding isomorphism 
restricts to the given one over $V$. 
\end{tm}
{\em Proof.}
Denote by $g$ the map $f_{|}: f^{-1} (V) \to V$. By cohomology and 
base change, $( Rf_{*} \rat_{X}[n] )_{|V} \simeq Rg_{*}  \rat_{f^{-1} 
(V) }[n]$. Similarly for $V'$.
Clearly, we have 
$( Rf_{*} \rat_{X}[n])_{|V}  \simeq \alpha^{*} 
[ (Rf_{*} \rat_{X} [n])_{|V'} ]$. 

\noindent
There is a   distinguished ``attaching'' triangle, see \ci{borel}, 
5.14:
$$
\beta_{*} \beta^{!}( Rf_{*} \rat_{X} [n]_{|V'} ) 
\stackrel{u}\lorw (Rf_{*} \rat_{X} [n])_{|V'} \stackrel{v}
\lorw R \alpha_{*} (Rf_{*} 
\rat_{X} [n]_{|V}) \stackrel{w  \; \; 
[1]}\lorw \beta_{*} \beta^{!}(Rf_{*} \rat_{X} [n]_{|V'})  [1].
$$
On the open set $V$ the complex $\beta_{*} \beta^{!}(Rf_{*} \rat_{X} [n])_{|V'} $ 
is isomorphic  to zero and 
the map $v$ restricts to  an isomorphism.
Recalling the notation  in $\S$\ref{bfars},
the long exact sequence of  cohomology sheaves 
is, stal-kwise along  the points of $S$, the long exact sequence for the 
cohomology of the pair $(f^{-1} (U ) , f^{-1} (U \setminus S')$.
In addition the map 
${\cal H}^{-d} (u)$ is identified, stalk-wise along the points of $S$, with the
map $\rho_{S,s}$.
The statement follows from Proposition \ref{rocom} which allows us to apply
proposition \ref{splits}.
\blacksquare

 \begin{rmk}
     {\rm 
    In the algebraic and K\"ahler case, the results \ci{bbd} and 
    \ci{sa}, coupled with Theorem \ref{fin}, 
 imply that the forms $\rho_{S}$ are non-degenerate for every $i \in I$;  see also
 \ci{cg}, Theorem 8.9.14. 
 To our knowledge these results have no implications as to the 
 sign of the intersection  forms. Surprisingly, in the projective case we can 
 determine that these forms are  definite; see $\S$\ref{sdtpc}.
 }
 \end{rmk}

 \medskip
 The following statement is known as a consequence of the 
characterizing properties  of  intersection cohomology sheaves; see
\ci{go-ma2}, $\S$6.2 and \ci{b-m}, $\S$1.8.
Recall that a map is {\em small} if all the inequalities in the 
definition of semismall are strict. Note that in this case there is 
only one relevant stratum.
 
\begin{cor}
\label{small}   
 Let  $f: X \to Y$ be a proper holomorphic small map from a connected complex 
 manifold.
 Then the Decomposition Theorem holds for $f$.
    \end{cor}
{\em Proof.} All the maps $\rho_{S}$ are identically zero and are 
isomorphisms.
\blacksquare

\subsection{Signature and Decomposition Theorem in the projective case}
\label{sdtpc}
In this section we use Theorem \ref{fai}, the previous inductive analysis 
and a Bertini-type argument to give a proof
of the following theorem

\begin{tm}
\label{dtp}
Let $f: X \to Y$ be a semismall map from a nonsingular complex
projective variety of dimension $n$ onto a  complex projective variety.
The Decomposition Theorem holds for $f$, i.e. there is a canonical isomorphism
 $$
  Rf_* \rat_{X} [n] \simeq \bigoplus_{ k \in I' }{
  IC_{\overline{Y_{k}}}   ({\cal  L}_{k} )} \simeq 
  \bigoplus_{ k \in I' }\bigoplus_{m=1}^{m_k}{
  IC_{\overline{Y_{k}}}   ({\cal  L}_{km} )}.
  $$
For every relevant stratum $S$  of dimension $d$ 
the associated intersection form is non-degenerate
and $(-1)^{\frac{n-d}{2}}$-definite.
\end{tm}  
{\em Proof.} By virtue of Proposition \ref{fin} we are reduced 
to checking that the intersection form
associated with a relevant d-dimensional stratum $S$  is non-degenerate 
and $(-1)^{\frac{n-d}{2}}$-definite.

\noindent
If $d=0$, then the conclusion follows from Proposition \ref{fai}.

\noindent
Let $d>0$. Let $A$ be a very ample divisor on $Y$. The line bundle $M:=f^*A$
is lef and generated by its global sections. By virtue of Proposition \ref{ber},
we can choose $d$ general sections $H_1, \ldots,H_d$ in the linear 
system $|A|$ such that their common 
zero locus $H$ 
has the property that $f^{-1}(H)$ is nonsingular of dimension 
$n-d$, $f^{-1}(H) \to H $ is semismall, $H$ meets $S$ at a non-empty
finite set of points $s_1, \ldots , s_r$ so that, for at least one index $1
\leq l \leq r$,  a small
neighborhood of a point $s_l$ in $H$ is a normal slice to $S$ at
$s_l$.
 By virtue 
of Theorem \ref{fai} the intersection form of $f^{-1}(s_l) \subseteq f^{-1}(H)$ 
has the required properties at the point $s_l$, and therefore
at every point $s \in S$. We conclude by applying Lemma \ref{itssymm}.
\blacksquare

\begin{rmk}
\label{elenco}
{\rm 
Theorem \ref{dtp} can be applied even when the
spaces are not complete, in the presence of a suitable completion
of the morphism: one for which the domain is completed to a projective
manifold,
the target to a projective variety and the map to a semismall one.
In general this may not be possible,
but it can be done in several instances,
e.g. the Springer resolution of the nilpotent cone of a complex semisimple 
Lie algebra, the Hilbert scheme of points on an algebraic surface
mapping on the corresponding symmetric product, isolated singularities
(see below), certain contraction of holomorphic symplectic varieties $\ldots$
}
\end{rmk}

\bigskip

\begin{cor}
 \label{isolatedsing}
Let $f: X \to Y$ be a birational   semismall map from a nonsingular 
quasi projective  complex variety of dimension $n$ onto a quasi projective 
complex variety $Y$ with isolated singularities. Assume that $f$ is 
an isomorphism outside the isolated  singularities.
The Decomposition Theorem holds for $f$. 
\end{cor}
{\em Proof.} We can reduce the statement to the complete projective 
case: see  Corollary \ref{issing} and Remark \ref{anynumber}.
\blacksquare



\bigskip
\noindent
Authors' addresses:

\smallskip
Mark Andrea A. de Cataldo,
Department of Mathematics,
SUNY at Stony Brook,
Stony Brook,  NY 11794, USA. \quad 
e-mail: {\em mde@math.sunysb.edu}

\smallskip
Luca Migliorini,
Dipartimento di Matematica, Universit\`a di Trento,
Via Sommarive, 14,
38050 Povo (Tn),  ITALY. \quad
e-mail: {\em luca@alpha.science.unitn.it}


\begin{thebibliography}{99}


\bibitem{bbd}{A.A. Beilinson, J.N. Bernstein, P. Deligne,
{\em Faisceaux pervers}, Ast\'erisque {\bf 100}, Paris, Soc. Math. Fr. 1982.}


\bibitem{borel} A. Borel et al., {\em Intersection Cohomology}, 
Progress in Mathematics Vol. {\bf 50}, Birkh\"auser, Boston Basel 
Stuttgart 1984.

\bibitem{b-m}{W. Borho, R. MacPherson, ``Partial resolutions of nilpotent 
varieties," Ast\'erisque {\bf 101-102} (1983), 23-74.}



\bibitem{cg}  N. Chriss, V. Ginzburg, Representation Theory and 
Complex Geometry, Birkh\'auser, Boston Basel 
Berlin 1997.

\bibitem{de-mi} M. de Cataldo, L. Migliorini, ``The Douady space
of a complex surface,"  Adv. in Math. {\bf 151}, 283-312 (2000).



\bibitem{dess} P. Deligne, ``Th\'eor\`eme de Lefschetz et crit\`eres de d\'eg\'en\'erescence
de suites spectrales,'' Publ.Math. IHES {\bf 35} (1969), 107-126.


\bibitem{ho2}  P. Deligne, ``Th\'eorie de Hodge, II,'' Publ.Math. IHES {\bf 40} (1971), 5-57.



\bibitem{weil2} P. Deligne, ``La conjecture de Weil, II,'' Publ.Math. IHES {\bf 52} (1980), 138-252.

\bibitem{e-v} H. Esnault, E. Viehweg, ``Vanishing and Non-Vanishing Theorems,"
Actes du Colloque de Th\'eorie de Hodge, 
Asterisque {\bf 179-180} 1989, 97-112.


\bibitem{fulton} W. Fulton, {\em Intersection Theory},
 Ergebnisse der Mathematik,
3.folge. Band 2, Springer-Verlag, Berlin Heidelberg 1984.


\bibitem{go-ma2} M. Goresky, R. MacPherson, ``Intersection homology II,''
Inv. Math. {\bf 71} (1983), 77-129. 


\bibitem{g-m} M. Goresky, R. MacPherson, {\em Stratified Morse Theory},
 Ergebnisse der Mathematik,
3.folge. Band 2, Springer-Verlag, Berlin Heidelberg 1988.



\bibitem{iv} B. Iversen, {\em Cohomology of Sheaves}, Universitext, 
Springer-Verlag, Berlin Heidelberg 1986.

\bibitem{laufer} H. Laufer, {\em Normal Two-Dimensinal Singularities},
Annals of Mathematics Studies {\bf 71}, Princeton University Press 1971.


\bibitem{looj} E. Looijenga, ``Cohomology
and Intersection Hohomology
of Algebraic Varieties,'' in {\em Complex Algebraic Geometry}, IAS/Park City 
Mathematics Series, Vol.3 J. Koll\'ar Editor, American Mathematical 
Society, 1997, 221-264.



\bibitem{mig} L. Migliorini, ``A smooth family of minimal surfaces of 
general type over a curve of genus at most one is trivial''
J.Algebraic Geometry {\bf 4}, 353-361 (1995)


\bibitem{sa}{M. Saito, ``Decomposition theorem for proper K\"ahler morphisms,"
Tohoku Math. J. (2) {\bf 42}, no. 2,  (1990),  127--147.}


\bibitem{wisn} J.A.Wi\'sniewski, ``Cohomological invariants of complex 
manifolds coming from extremal rays,''
Asian J.Math. Vol.2, No.{\bf 2}(1998), 289-302
\end{thebibliography}
\end{document}